\documentclass[12pt]{article}
\usepackage{natbib,graphicx,setspace,lscape,longtable}
\usepackage{natbib,epsfig,graphicx,float}
\usepackage{amsmath,amsthm,amssymb,color}
\usepackage[nohead]{geometry}
\usepackage[bottom]{footmisc}
\usepackage{indentfirst}
\usepackage{endnotes}
\usepackage{verbatim}
\usepackage{threeparttable}
\usepackage{comment}
\usepackage{algorithm}
\usepackage{titlesec}
\usepackage{multibib}
\usepackage{algorithmic}
\usepackage{multirow}
\usepackage{dsfont}
\usepackage{titling}
\usepackage{afterpage}
\usepackage{hyperref}
\RequirePackage[mathlines, displaymath]{lineno}
\usepackage{multirow}
\usepackage{chngcntr}
\usepackage[titletoc]{appendix}
\bibpunct{(}{)}{;}{a}{,}{,}
\floatstyle{ruled}
\newfloat{algorithm}{tbhp}{loa}
\floatname{algorithm}{Algorithm}

\numberwithin{equation}{section}
\theoremstyle{plain}

\theoremstyle{definition}

\newtheorem{theorem}{\bf Theorem}
\newtheorem{lemma}{\bf Lemma}

\newtheorem{assumption}{\bf Assumption}

\def\var{\mathrm{var}}
\def\tsum{\textstyle{\sum}}

\def\two{I\!I}
\def\three{I\!I\!I}
\def\four{I\!V}

\newcommand{\hf}{\widehat f}
\newcommand{\hB}{\widehat B}
\newcommand{\hF}{\widehat F}
\newcommand{\hLam}{\widehat \Lambda}
\newcommand{\hSigma}{\widehat \Sigma}
\newcommand{\hg}{\widehat g}
\newcommand{\hphi}{\widehat\phi}
\def\inv{^{-1}}
\def\half{^{1/2}}
\def\nano{\scriptscriptstyle}
\newcommand{\indep}{\rotatebox[origin=c]{90}{$\models$}}
\newcommand\lo[1]{_{\nano #1}}
\def\cstrue{{\mathcal S}\lo {y | f}}
\def\csest{{\mathcal S}\lo {y | \hf}}

\def\bop{O_P}
\def\sop{o_P}
\def\tsum{\textstyle\sum}

\allowdisplaybreaks[4]
\begin{document}
\title{Inverse Moment Methods for Sufficient Forecasting using High-Dimensional Predictors}
\author{Wei Luo$^1$, Lingzhou Xue$^2$, Jiawei Yao$^3$ and Xiufan Yu$^2$ \\ $^1$Zhejiang University, $^2$Pennsylvania State University and $^3$Princeton University
}
\date{}
\maketitle{}

\begin{abstract}
We consider forecasting a single time series using a large number of predictors in the presence of a possible nonlinear forecast function. {Assuming that the predictors affect the response through the latent factors, we propose to first conduct factor analysis and then apply sufficient dimension reduction on the estimated factors, to derive the reduced data for subsequent forecasting.} Using directional regression and the inverse third-moment method in the stage of sufficient dimension reduction, the proposed methods can capture the non-monotone effect of factors on the response. We also allow a diverging number of factors and only impose general regularity conditions on the distribution of factors, avoiding the undesired time reversibility of the factors by the latter. These make the proposed methods fundamentally more applicable than the sufficient forecasting method in \cite{fxy-2015}. The proposed methods are demonstrated in both simulation studies and an empirical study of forecasting monthly macroeconomic data from 1959 to 2016. Also, our theory contributes to the literature of sufficient dimension reduction, as it includes an invariance result, a path to perform sufficient dimension reduction under the high-dimensional setting without assuming sparsity, and the corresponding order-determination procedure.
\end{abstract}


\pagestyle{plain}

\section{Introduction}

Forecasting using high-dimensional predictors is an increasingly important research topic in statistics, biostatistics, macroeconomics and finance. A large body of literature has contributed to forecasting in a data rich environment, with various applications such as the forecasts of market prices, dividends and bond risks \citep{sharpe-1964,lintner-1965,Ludvigson-Ng-2009}, macroeconomic outputs \citep{stock-watson-1989,bernanke-etal-2005}, macroeconomic uncertainty and fluctuations \citep{Ludvigson-Ng-2007,jurado-2015}, and clinical outcomes based on massive genetic, genomic and imaging measurements. Motivated by principal component regression, the pioneering papers by \cite{stock-watson-2002,stock-watson-2002b} systematically introduced the forecasting procedure using factor models, which has played an important role in macroeconomic analysis. Recently, \cite{fxy-2015} extended \cite{stock-watson-2002,stock-watson-2002b} to allow for a nonlinear forecast function and multiple nonadditive forecasting indices. Following \cite{fxy-2015}, we consider the following factor model with a target variable $y_{t+1}$ that we aim to forecast:
\begin{eqnarray}
y_{t+1} &=& g(\phi_1'f_t,\cdots, \phi_L'f_t,\epsilon_{t+1}),  \label{eq1.1}\\
x_{it} &=& b_i'f_t+u_{it},\quad 1\leq i\leq p,\ 1\leq t\leq T, \label{eq1.2}
\end{eqnarray}
where $x_{it}$ is the $i$-th high-dimensional predictor observed at time $t$, $b_i$ is a $K\times 1$ vector of factor loadings, $f_t$ is a $K\times 1$ vector of common factors driving both predictor and response, $g(\cdot)$ is an unknown forecast function that is possibly nonadditive and nonseperable, $u_{it}$ is an idiosyncratic error, and $\epsilon_{t+1}$ is an independent stochastic error. Here, $\phi_1,\ldots,\phi_L$, $b_1,\ldots,b_p$ and $f_1,\ldots,f_T$ are unobserved vectors. Model (\ref{eq1.1}) equivalently assumes

\begin{equation}\label{eq:sdr}
y_{t+1} \indep f_t \, | \, (\phi_1,\ldots, \phi_L)' f_t.
\end{equation}

The linear space spanned by $\phi_1,\ldots, \phi_L$, denoted by $\cstrue$, is the parameter of interest that is identifiable and known as the central subspace \citep{cook1998}. \cite{fxy-2015} introduced the sufficient forecasting scheme to use factor analysis in model (\ref{eq1.2}) to estimate $f_t$, and apply the sliced inverse regression \citep{li-1991} in model (\ref{eq1.1}) with the estimated factors as the predictor. Such a combination provides a promising forecasting technique that not only extracts the underlying commonality of the high-dimensional predictor but also models the complex dependence between the predictor and the forecast target. It allows the dimension of the predictor to diverge and even become much larger than the number of observations.%, which is intrinsically appealing to solving high-dimensional forecasting problems.

The consistency result of \cite{fxy-2015} is not granted as it may appear. If we replace the true factors $f_t$ with a {consistent} estimate $\hf_t$ in (\ref{eq:sdr}) and define the central subspace $\csest$ similarly, then $\csest$ may differ with $\cstrue$ drastically in general. Thus, the naive method by applying existing dimension reduction methods to the estimated factors $\hf_t$'s may not necessarily lead to the consistent estimation of $S_{y|f}$, {even if it consistently estimates $\csest$.} \cite{fxy-2015} effectively {addressed} this issue by developing an important invariance result between $E(f_t | y_{t+1})$ and $E(\hf_t|y_{t+1})$. See Proposition 2.1 and Equation (2.9) of \cite{fxy-2015}. This invariance result provides an essential foundation for using the sliced inverse regression under Models \eqref{eq1.1}--\eqref{eq1.2}.

{Nonetheless, the applicability of \cite{fxy-2015} is restricted by the requirements that the number of factors $K$ must be fixed as $p$ and $T$ grow, and, for each set of factors, a linearity condition (see (B1) below) must hold. In particular, as $\cstrue$ is unknown, the linearity condition is commonly strengthened}  {to equivalently require an elliptically distributed $f_t$, which causes the undesired time reversibility \citep{xia-etal-2002}. In addition, the consistency of \cite{fxy-2015} and \cite{yu2020} hinges on an exhaustive estimation of $\cstrue$, (i.e. detecting all the directions), for which $\phi_1'\Sigma_{f|y} \phi_1, \ldots, \phi_L'\Sigma_{f|y} \phi_L$ must be positive (See their Assumption (A2)). This condition is violated, i.e. $\phi' \Sigma_{f|y} \phi$ being zero for some $\phi \in \cstrue$, if $\phi' f_t | y_{t+1}$ has a symmetric distribution, which occurs when the forecast target was investigated using squared factors \citep{bai-ng-2008,Ludvigson-Ng-2007}. These limitations motivate us to construct more powerful forecasting methods based on \cite{fxy-2015}'s work.}

{In this paper, we propose to use factor analysis and sufficient dimension reduction sequentially for sufficient forecasting, with second- or higher-order inverse moment methods being the working sufficient dimension reduction method. In the main text, we focus on a commonly used second-order inverse moment method called directional regression \citep{li-wang-2007}, and defer the development with the third-order inverse moment method to the online supplement. Based on models (\ref{eq1.1}) and (\ref{eq1.2}), the proposed method includes the following steps:}

\noindent
Step 1. Estimate the factor loadings $B$ and the factors $f_t$ in Model (\ref{eq1.2}).

\noindent
Step 2. Use the estimates $\widehat B$ and $\hf_t$ in directional regression to estimate $\cstrue$.

\noindent
Step 3. Use the nonparametric methods \citep{fan-1996, matzkin-2003, yu2020} to estimate $g(\cdot)$ in Model (\ref{eq1.1}) and forecast $y_{t+1}$, based on the estimate of $(\phi_1'f_t,\ldots, \phi_L'f_t)$.

By studying both $E(f_t|y_{t+1})$ and $E(f_tf_t'|y_{t+1})$ in Step 2, we explore the full power of the factor space. To this end, we first provide an important invariance result (i.e. Lemma \ref{lemma:invariance}) for {directional regression}. With the help of this invariance result, we do not require the coincidence or closeness of two central subspaces $S_{y| f}$ and $\csest$, so the proposed method can be applied to more general data, such as non-normally distributed factors.

Our work extends the method, theory and applicability of the forecasting using factor models. Compared with \cite{fxy-2015}, we relax the linearity condition to the general moment conditions on $f_t$. From the discussion above, the proposed method does not require time reversibility of the factors, so it can be applied to the generalized forecasting model
\begin{equation}\label{eq: model lag variable}
y_{t+1}=g(\phi_1'f_t+\psi_1'\omega_t, \ldots, \phi_L'f_t+\psi_L'\omega_t,\epsilon_{t+1})
\end{equation}
where $\omega_t$ is an $m\times 1$ vector of the observed variables (e.g. lags of $y_{t+1}$).
In addition, by using the higher-order inverse moments, the proposed method requires weaker condition than \cite{fxy-2015} and \cite{yu2020} for exhaustive estimation of $\cstrue$. In particular, it can detect non-monotone effect of the factors on the response. Furthermore, we allow the number of underlying factors $K$ to diverge as $p,T\rightarrow \infty$. By \cite{lam-yao-2012,li-etal-2013} and \cite{jurado-2015}, our method will deliver a more powerful forecast than \cite{stock-watson-2002,stock-watson-2002b} and \cite{fxy-2015}.

Using the directional regression as an illustration, the proposed method also provides a novel framework of performing sufficient dimension reduction with large panel data under the high-dimensional setting, without the commonly-adopted sparsity assumption but with the assumption that the predictor affects the response only through the latent factors. The original direction regression \citep{li-wang-2007} can only deal with independently and identically distributed data under the low-dimensional setting. This enhances the applicability of model-free dimension reduction for high-dimensional data, when the sparsity assumption is not suitable.

%We have also developed the methodology with the inverse third-moment methods used in Step 2 above, which requires even weaker condition for exhaustive estimation of $\cstrue$. To save space, this is deferred to the supplement.

The consistency of the proposed method hinges on the consistency of both factor analysis and directional regression based on the estimated factors, which we study next. For ease of presentation, we assume that both the number of factors $K$ and the dimension $L$ of $\cstrue$ known a priori. This does not affect the asymptotic development of the resulting estimator, as long as $K$ and $L$ can be consistently estimated; see the supplement for details. The consistent estimation of $K$ and $L$ is deferred to \S 5. {Throughout the article, we assume $L$ to be fixed as $K$ diverges.}

%The rest of this paper is organized as follows. We first study the sufficient forecasting with directional regression in Section 2, including factor analysis in Subsection 2.1, an invariance result in Subsection 2.2, the details of implementation in Subsection 2.3, asymptotic results in Subsection 2.4, relaxation of the linearity condition in Subsection 2.5, and model selection in Subsection 2.6. In Section 3, we study the inverse third-moment method in the sufficient forecasting. Section 4 is devoted to the simulation studies and a real data example that illustrates the power of the proposed method. We leave technical proofs to Section 5 and the supplement. Section 6 includes a few concluding remarks.

%\section{Sufficient forecasting with directional regression}

\section{Consistency of factor analysis}

To make forecast, we need to estimate the factor loadings $B$ and the error covariance matrix $\Sigma_u$. Consider the following constrained least squares problem:
\begin{eqnarray}\label{eq:est B and F}
(\hB_K,\hF_K) =\text{arg}\min\limits_{(B,F)}&\|X-BF'\|^2_F,
\text{ subject to } T^{-1}F'F=I_K, \,\, B'B \text{ is diagonal,}
\end{eqnarray}
where $X=(x_1,\cdots,x_T), F'=(f_1,\cdots,f_T)$, and $\|\cdot\|_F$ denotes the Frobenius norm of a matrix. The constraints $T^{-1}F'F=I_K$ and that $B'B$ is diagonal address the issue of identifiability during the minimization. As these conditions can always be satisfied for any $BF'$ after appropriate matrix operations on $B$ and $F$, they impose no additional restrictions on the factor model (\ref{eq1.2}). It is known that the minimizers $\hF_K$ and $\hB_K$ of (\ref{eq:est B and F}) are such that the columns of $\hF_K/\sqrt{T}$ are the eigenvectors corresponding to the $K$ largest eigenvalues of the $T\times T$ matrix $X'X$ and $\hB_K=T^{-1}X\hF_K$. To simplify notation, let $\hB=\hB_K$ and $\hF=\hF_K$.

As both the dimension $p$ of the predictor $x_t$ and the number of factors $K$ are diverging, it is necessary to regulate the magnitude of the factor loadings $B$ and the idiosyncratic error $u_t$, so that the latter is negligible with respect to the former. We should also regulate the stationarity of the time series. In this paper, we adopt the following assumptions. For simplicity in notation, we let $U=(u_{it})_{p \times T}$, $B = (b_1,\ldots, b_p)'$, and $\|B\|_{\max}$ be the maximum of the absolute values of all the entries in $B$. Let $\mathcal{F}_{\infty}^0$ and $\mathcal{F}_T^{\infty}$ denote the $\sigma-$algebras generated by $\{(f_t,u_t,\epsilon_{t+1}): t\leq 0\}$ and $\{(f_t,u_t,\epsilon_{t+1}): t\geq T\}$ respectively. Let $\alpha(T)=\sup_{A\in\mathcal{F}_{\infty}^0, B\in\mathcal{F}_T^{\infty}} |P(A)P(B)-P(AB)|$.

\begin{assumption}[Factors and Loadings] \label{assume:factors and loadings} \empty \hfill
\newline
(1) There exists $b > 0$ such that $\|b_i\|\leq b$ for $i=1,\ldots,p$
%$\sup_{p \in \mathbb{N}} \|\bB\|_{\max} \leq b$,
and there exist two positive constants $c_1$ and $c_2$ such that $
c_1< p\inv \lambda_{\min}(B'B)< p\inv \lambda_{\max}(B'B)<c_2$;\\
(2) Identification: $T^{-1}F'F =I_K$, and $B'B$ is a diagonal matrix with distinct entries.
%(3) Linearity: $E(\bb'\bff_t|\bphi_1'\bff_t,\cdots,\bphi_L'\bff_t)$ is a linear function of $\bphi_1'\bff_t,\cdots,\bphi_L'\bff_t$ for  any $\bb\in\R^p$.
\end{assumption}

\begin{assumption}[Data Generating Process] \label{assume:data generating}
$\{f_t\}_{t\geq 1}$, $\{u_t\}_{t\geq 1}$ and $\{\epsilon_{t+1}\}_{t\geq 1}$ are three independent groups, and they are strictly stationary. $\{K^{-2}E\|f_t\|^4:  K \in \mathbb{N}\}$ and $\{ K^{-1}E(\|f_t\|^2|y_{t+1}):  K\in \mathbb{N}\}$ are bounded sequences, and $\alpha(T)<c\rho^T$ for $T\in\mathbb{Z}^+$ and some $\rho\in(0,1)$.
\end{assumption}
\begin{assumption}[Residuals and Dependence] \label{assume: dependence}
There is a constant $M>0$ %that does not depend on $p$ or $T$
such that
(1) %$E(\bu_t)=\bzero$ and
$E|u_{it}|^8\leq M$;
(2) $\|\Sigma_u\|_1\leq M$; %,  and for every $i,j,t,s>0$, $(pT)^{-1}\sum_{i,j,t,s}|E(u_{it}u_{js})|\leq M$; \
%(2) $\bSigma_u=(\sigma_{ij})_{p\times p}$ satisfies (a) $p^{-1}\sum_{i=1}^p|\sigma_{ii}|\leq M$, and (b) $p^{-1}\sum_{i,j}|\sigma_{i,j}|\leq M$. \\
%$E(u_{it}u_{jt})=\sigma_{ij,t}$. $|\sigma_{ij,t}|\leq \bar\sigma_{ij}$ for some $\bar\sigma_{ij}>0$, and\\
%\indent(a) $|p^{-1}\sum_i\sigma_{ii,ss}|\leq M$, (b) $p^{-1}\sum_{i,j}\bar\sigma_{ij}\leq M$, (c) $(pT)^{-1}\sum_{i,j,s,t}|\sigma_{ij,st}|\leq M$.  \\
(3) For every $(t,s)$, $E|p^{-1/2}( u_s'u_t-E(u_s' u_t) )|^4\leq M$;
%(4) Weak dependence between factors and idiosyncratic errors
% $$E(\frac{1}{p}\sum_{i=1}^p\|\frac{1}{\sqrt{TK}}\sum_{t=1}^T\bff_tu_{it}\|^2)\leq M$$
(4) $U=LER$ where $L\in\mathbb{R}^{p\times p}$ and $R\in\mathbb{R}^{T\times T}$ are non-random positive definite matrices and $E=(e_{it})_{p\times T}$ includes independent elements with $E(e_{ti})=0$ and $E|e_{it}|^7\leq M$.
\end{assumption}

Assumptions \ref{assume:factors and loadings} and \ref{assume: dependence} ensure that signals dominate errors in the population level as $p$ grows. Assumptions \ref{assume:factors and loadings} regulates the signal strength of factors contained in the predictor through the convergence rate of estimated factor loadings, and Assumption \ref{assume: dependence} regulates the idiosyncratic errors. Assumption \ref{assume: dependence}(4) regulates weak autocorrelation and cross-sectional correlation as in \cite{li-etal-2013}. Assumption \ref{assume:data generating} imposes independence between factors and idiosyncratic errors as in \cite{lam-yao-2012}. Assumption \ref{assume:data generating} implies that the observations are only weakly dependent, so that the estimation accuracy grows with $T$. Assumption \ref{assume:data generating} and Assumption \ref{assume: dependence}(2) imply that for every $i,j,t,s>0$, $\max_{t\le T}p^{-1}\sum_{i,j}|E(u_{it}u_{jt})|=O(1)$ and $(pT)^{-1}\sum_{i,j,t,s}|E(u_{it}u_{js})|=O(1)$ (See Lemma 6 of \cite{fan-2013}).

Under these assumptions, we have the following consistency result for estimating the factor loadings. Instead of the Frobenius norm used in (\ref{eq:est B and F}), we use the spectral norm to measure the magnitude of a matrix, defined as $\|A\| = \lambda_{\max}^{1/2}(A'A) $, the square root of the largest eigenvalue of the symmetric matrix $A'A$, for any matrix $A$.

\begin{theorem}\label{thm:order of Bhat and Lambda}
Let $\Lambda_b = (B'B)\inv B'$ and $\hLam_b = (\hB'\hB)\inv \hB'$. Given $K = o(\min\{p^{1/3},T\})$ and Assumptions \ref{assume:factors and loadings}, \ref{assume:data generating} and \ref{assume: dependence}(1)-(3), we have
\begin{itemize}
\item[1)] $\|\hB-B\|=O_p(p^{1/2}(K^{3/2} p^{-1/2} + K^{1/2} T^{-1/2}))$,
\item[2)] $\|\hLam_b-\Lambda_b\|=O_p(p^{-1/2}(K^{3/2} p^{-1/2} + K^{1/2} T^{-1/2}))$.
\end{itemize}
\end{theorem}

Theorem \ref{thm:order of Bhat and Lambda} extends the existing consistency result for estimating the factor loadings \citep{lyb-2011, fan-2013, fxy-2015} by pinpointing the effect of diverging $K$. Because the dimension $p$ of factor loadings $B$ is diverging, the estimation error $\hB - B$ accumulates as $p$ grows. For a $p$-dimensional vector whose entries are constantly one, its spectral norm is $p\half$, which diverges to infinity. Thus, we should treat $p\half$ as the unit magnitude of the spectral norm of matrices with $p$ rows, in which sense the statement 1) of Theorem \ref{thm:order of Bhat and Lambda} justifies the estimation consistency of the factor loadings $B$. As the error term $u_t$ shrinks as $p$ grows under Assumption \ref{assume: dependence}, the convergence rate of the factor loading estimation largely depends on $p$ - a higher dimensional predictor means a more accurate estimation. The convergence rate in this theorem can be further improved if we impose stronger assumptions on the negligibility of the error terms in the factor model (\ref{eq1.2}).

Given $\hB$, it is easy to see $\widehat f_t = \hLam_b B f_t + \hLam_b u_t$. Thus, together with the negligibility of the error term $u_t$, the consistency of $\hB$ and $\hLam_b$ indicates the closeness between the true factors $f_t$ and the estimated factors $\widehat f_t$, of which the latter will be used in the subsequent sufficient dimension reduction. The error covariance matrix $\Sigma_u$ can be estimated by thresholding the sample covariance matrix of the estimated residual $x_t - \hB \widehat f_t$, denoted by $\hSigma_{u}=(\hat\sigma^{u}_{ij})_{p\times p}$, as in \cite{cai-liu-2011}, \cite{xue-ma-zou-2012} and \cite{fan-2013,fxz-2015}.
%\begin{equation}\label{eq:est erro variance-1}
%\hSigma^{\mathcal T_1}_{u}=\left(\hat\sigma^{u}_{ij}I_{\{i=j\}}+s_{\tau_{ij}}(\hat\sigma^{u}_{ij})I_{\{i\neq j\}}\right)_{p\times p}
%\end{equation}
%where $s_{\tau_{ij}}(\cdot)$ is the soft-thresholding rule and ${\tau_{ij}}=C\omega_{\mathcal T} \hat\theta_{ij}\half >0$ is the adaptive threshold for each entry. Alternatively, we may employ \cite{xue-ma-zou-2012} to solve the positive definite sparse estimator from
%\begin{equation}\label{eq:est erro variance-2}
%\hSigma^{\mathcal T_2}_{u}=\arg\min_{\boldsymbol{\Sigma} \succeq \epsilon \bI} \ \frac{1}2\|\bSigma-\hat\bSigma_{u}\|_F^2 \ +\sum_{(i,j):i\neq j}\tau_{ij}|\sigma_{ij}|,
%\end{equation}
%where $A \succeq B$ means the matrix $A - B$ is positive semi-definite for square matrices $A$ and $B$.
%For simplicity, denote the sparse error covariance estimator by $\hSigma^{\mathcal T}_{u}$.

\section{Directional regression based on an invariance result}

\subsection{An invariance result}

Had the true factors $f_t$ been observed, directional regression would estimate the central subspace $\cstrue$ as the column space of
\begin{equation}\label{eq:kernel-dr} M_{dr} = E\{2 \var(f_t) - E[(f_t - g_s) (f_t - g_s)' | y_{t+1}, \eta_{s+1}]\}^ 2, \end{equation}
where $(g_s, \eta_{s+1})$ is a hypothetical independent copy of $(f_t, y_{t+1})$. The term $\var(f_t)$ can be replaced with the identity matrix as in \cite{li-wang-2007}, but we keep it in this form for the convenience in the theoretical work developed later.  For the resulting directions being included in $\cstrue$, $f_t$ needs to satisfy the linearity condition and the constant variance condition; that is,

\noindent (B1) $E(b'f_t|\phi_1'f_t,\cdots,\phi_L'f_t)$ is a linear function of $(\phi_1' f_t,\ldots,\phi_L' f_t)$ for any $b \in \mathbb {R}^K$.

\noindent (B2) $\var(f_t \mid \phi_1'f_t,\ldots, \phi_L'f_t)$ is degenerate.

\noindent Since $\cstrue$ is unknown, (B1) and (B2) are commonly strengthened such that they are satisfied for basis matrices of any $L$-dimensional subspace of $\mathbb {R}^K$. The strengthened conditions equivalently require the factors to be jointly normally distributed. To assess these conditions, one can treat $f_t$ as the response and $(\phi_1' f_t,\ldots,\phi_L'f_t)$ as the predictor in regression, then (B1) is the linearity assumption on the regression function and (B2) is the homoscedasticity assumption on the error term. In this sense, we follow the convention in the literature of regression to treat (B2) less worrisome than (B1) in practice. {We tentatively assume (B1) and relax it in \S 4.} %In this paper, we would assume that (B2) holds.

Under general conditions, the column space of $M_{dr}$ is $L$-dimensional, which, together with the linearity condition (B1) and the constant variance condition (B2), means the exhaustive recovery of $\cstrue$. These conditions are proposed in \cite{li-wang-2007} and reviewed in the supplement. They are weaker than those required for the exhaustiveness of sliced inverse regression, as more information about $f_t | y_{t+1}$, i.e. the second moment, is used. We assume these conditions throughout the paper, including \S 4 where (B1) is violated.

To pinpoint the effect of using the estimated factors in directional regression, we next propose an invariance result for $M_{dr}$. As mentioned in \S 1, a similar invariance result for sliced inverse regression can be found in \cite*{fxy-2015} where only the inverse first moment is involved (see their equation (2.6)). To simplify the discussion, in the rest of the subsection we assume an oracle scenario where $B$ is known a priori, which gives
\begin{equation}\label{eq:hat f w known B}
\hf_t = f_t + u^*_t,
\end{equation}
where $u^*_t = \Lambda_b u_t$ is independent of $f_t$. Let $u_s^*$ be an independent copy of $u_t^*$ in (\ref{eq:hat f w known B}) and let $\hg_s = g_s + u_s^*$. Since $B$ is known, $\widehat g_s$ is an independent copy of $\widehat f_t$.

\begin{lemma}\label{lemma:invariance}
(The invariance result) Under model (\ref{eq1.2}), $M_{dr}$ defined in (\ref{eq:kernel-dr}) is invariant if the true factors $f_t$ and $g_s$ are replaced with the estimated factors $\hf_t$ and $\hg_s$. %and $I_K + (\bB'\bB)^{-1} \bB' \Sigma_u \bB (\bB' \bB)^{-1}$, respectively.
\end{lemma}

Using the estimated factors, one would naturally treat $S\lo {y|\hf}$ as the working parameter in the stage of sufficient dimension reduction. However, as no distributional assumptions are imposed on $u^*_t$, both (B1) and (B2) can be violated for $\hf_t$, which causes inconsistency of directional regression for recovering $\csest$. In addition, $\csest$ itself may deviate from the parameter of interest $\cstrue$, as the identity between the two essentially requires the normality of both $f_t$ and $u^*_t$ \citep{li-yin-2007}. The invariance result provides the key to address these issues; that is, we can bypass $\csest$ and directly estimate $\cstrue$ using the estimated factors, as if the true factors were used. As $\var (\hf_t)$ is no longer the identity matrix, $M_{dr}$ adopted here modifies its original form in \cite{li-wang-2007}. This modification is crucial as it averages out the effect of the estimation error $u_t^*$. It also means that the column space of the working $M_{dr}$ does differ from $\csest$.

%By Theorem \ref{thm:order of Bhat and Lambda}, the factor loadings $\bB$ are consistently estimated by $\hB$, so it is plausible that under suitable moment conditions, the leading eigenvectors of the resulting matrix span a consistent estimator of the central subspace $\cstrue$. %The corresponding convergence rate will be given in Subsection 2.4 under Assumptions 2.1 - 2.3.

\subsection{Consistency of directional regression}

In reality, the hypothetical independent copies $(g_s, \eta_{s+1})$ and $(f_t, y_{t+1})$ do not exist in the observed data, so we expand
(\ref{eq:kernel-dr}) and estimate an equivalent form of $M_{dr}$,
\begin{eqnarray}\label{eq:kernel-dr-equiv}
M_{dr} &=&  2E\{[\var(f_t) - E (f_tf_t' |y_{t+1})]^2\} + 2E^ 2 [E(f_t|y_{t+1})E'(f_t |y_{t+1})] \nonumber \\
&& + 2E[E'(f_t |y_{t+1})E(f_t|y_{t+1})]\cdot E[E(f_t|y_{t+1})E'(f_t |y_{t+1})].
\end{eqnarray}
By Lemma \ref{lemma:invariance}, we can replace $f_t$ with $\hf_t$, in which $B$ is replaced with $\widehat B$. For the ease of estimation, in the literature {of sufficient dimension reduction}, it has been a common practice to employ the slicing technique; that is, we partition the sample of $y_{t+1}$ into $H$ slices with equal sample proportion. In the population level, it corresponds to partitioning the support of $y_{t+1}$ into $H$ slices with equal probability, and using the corresponding indicator, denoted by $y_{t+1}^ D$, as the new working response variable. %Equivalently, we estimate the parameter of interest, the central subspace $\cstrue$, via the intermediate parameter, the central subspace $\csest$.

Because the slice indicator $y_{t+1}^D$ is a measurable function of the original response $y_{t+1}$, $f_t$ must affect $y_{t+1}^D$ through $y_{t+1}$. Thus, the working central subspace ${\mathcal S}_{y^D | f}$ is always a subspace of the central subspace of interest $\cstrue$. The two spaces further coincide for large $H$. Because the dimension $L$ of $\cstrue$ is fixed as $K$ grows, without loss of generality, we fix $H$ as $K$ grows and assume the identity between ${\mathcal S}_{y^D|f}$ and $\cstrue$. Such identity is conformed by an omitted simulation study that shows the robustness of the proposed method to the choice of $H$, for a reasonable range of $H$, e.g. from three to ten. The same phenomenon has also been commonly observed in the literature \citep{li-1991, li-wang-2007}.

Using $y_{t+1}^D$, the inverse moments $E (\hf_t |y_{t+1})$ and $E (\hf_t\hf_t' |y_{t+1})$ in ${M}_{dr}$ become the marginal moments of $\hf_t$ within each slice, and can be estimated by the usual sample moments. Hence, the slicing technique simplifies the estimation. In detail, we have

\medskip

\noindent Implementation of Step 2. Let $y_{(0)/H} = -\infty$, and, for $i=1,\ldots, H$, let $y_ {(i)/H}$ be the $(i/H)$th quantile of $\{y_1,\ldots, y_T\}$. Let $y_{t+1}^D = i$ if $y_{t+1} \in ( y_{(i)/H}, y_{(i+1)/H} ]$. Estimate $E(\hf_t | y_{t+1}^D = i)$ by $\tsum_{t=1}^T \hf_{t} I(y_{t+1}^D = i)/ (T/H)$ and $E(\hf_t \hf_t' | y_{t+1}^D =i)$ by $\tsum_{t=1}^T \hf_{t}\hf_{t}' I(y_{t+1}^D = i)/ (T/H)$. Estimate $\var (\hf_t)$ by $I_K$. Plugging these into (\ref{eq:kernel-dr-equiv}) to derive $\widehat M_{dr}$. Estimate $\cstrue$ by the space spanned by $(\widehat \phi_1,\ldots, \widehat \phi_L)$, the leading $L$ eigenvectors of $\widehat M_{dr}$.

\medskip

\noindent To estimate $\var (\widehat f_t)$ in (\ref{eq:kernel-dr-equiv}), one can alternatively use $I_K + \widehat \Sigma_{u^*}$ by the restriction $\var (f_t) = I_K$, where $\widehat \Sigma_{u^*}$ is the thresholding covariance estimator. An omitted simulation study shows that the resulting estimator of $M_{dr}$ performs similarly.

\begin{theorem}\label{thm:dr-consist}
Suppose $K = o(\min(p^{1/3}, T^{1/2}))$. Under Assumptions \ref{assume:factors and loadings}, \ref{assume:data generating} and \ref{assume: dependence}(1)-(3), the linearity condition (B1), and the constant variance condition (B2),  $(\widehat \phi_1,\ldots,\widehat \phi_L)$ span a consistent estimator of $\cstrue$ in the sense that
\begin{equation*}
\| (\widehat \phi_1,\ldots,\widehat \phi_L) (\widehat \phi_1,\ldots,\widehat \phi_L)' -  (\phi_1,\ldots, \phi_L)  ( \phi_1,\ldots, \phi_L)'\|_F = \bop (K^{3/2}p^{-1/2} + KT^{-1/2}).
\end{equation*}
\end{theorem}

In connection with Theorem \ref{thm:order of Bhat and Lambda}, this theorem justifies that the estimation error of $\cstrue$ comes from two parts. The first part, which is of order $\bop (K^{3/2}p^{-1/2})$, is inherited from factor analysis. This part represents the price we pay for estimating the factor loadings $B$, and it depends on the dimension $p$ of the original predictor. By contrast, the second part, which is of order $\bop (K T^{-1/2})$, does not depend on $p$ and is newly generated in the sufficient dimension reduction stage. From the proof of Theorem \ref{thm:dr-consist} (see the supplement), it represents the price we pay for estimating the unknown inverse second moment involved in the kernel matrix. Therefore, this part would persist even if no error were generated in factor analysis.

\section{Relaxing the linearity condition}

As mentioned in \S 3, (B1) can be regarded as a parametric assumption and can be violated in real applications. For example, this occurs when one incorporates the lag variables of $y_{t+1}$ in forecasting and consider Model (\ref{eq: model lag variable}). In this section, we address this issue in two ways: first, we justify the consistency of the proposed method without (B1) under the setting that the number of factors $K$ must diverge; second, we weaken (B1) and generalize the proposed method accordingly following the spirit of \cite{dong2010} under the setting that $K$ is fixed. %These together enhance the applicability of the proposed method.

When (B1) is violated, Theorem \ref{thm:dr-consist} still holds if we treat $(\phi_1,\ldots, \phi_L)$ as the $L$ leading eigenvectors of $M_{dr}$. Thus, the consistency of the proposed methodology depends on the closeness between the column space of $M_{dr}$ and the central subspace $\cstrue$, which hinges on the approximation of (B1). Fortunately, the latter has been justified in \cite{hall-li-1993} for all large $K$.

\begin{theorem}\label{cor:dr-consist}
Suppose $K \rightarrow \infty$ and $K = o(\min(p^{1/3}, T^{1/2}))$. Under Assumptions \ref{assume:factors and loadings}, \ref{assume:data generating}, \ref{assume: dependence}(1)-(3), the constant variance condition (B2), and other regularity conditions (see the supplement), $\widehat \phi_1,\ldots,\widehat \phi_L$ span a consistent estimator of $\cstrue$ in the sense that
\begin{equation*}
\| (\widehat \phi_1,\ldots,\widehat \phi_L) (\widehat \phi_1,\ldots,\widehat \phi_L)' -  (\phi_1,\ldots, \phi_L)  ( \phi_1,\ldots, \phi_L)'\|_F = \sop (1).
\end{equation*}
\end{theorem}

In the literature, \citeauthor{hall-li-1993}'s result on the approximation of (B1) was used heuristically to support the effectiveness of inverse moment methods when (B1) is violated; see, for example, \cite{cook-1991} and \cite{li-wang-2007}. As we are aware of, this is the first attempt to rigorously build the consistency of inverse moment methods using \citeauthor{hall-li-1993}'s result.

When $K$ is small and the factors clearly violate (B1), the approximation result in \cite{hall-li-1993} no longer applies. In this case, we treat $K$ as fixed, and relax (B1) to

\noindent(B1') $E(f_t | \phi_ 1' f_t,\ldots, \phi_ L' f_t)$ is a linear combination of $\{h_i (\phi_1' f_t,\ldots, \phi_L' f_t): i=1,\ldots, q\}$.

\noindent
One can set the basis functions in (B1') to be power functions, trigonometric functions, etc. In addition to (B1'), we require the constant variance condition (B2), which, as mentioned in \S 1, is quite mild. These conditions closely resemble those in \cite{dong2010}. We generalize directional regression from the eigen-decomposition of $M_{dr}$ to minimizing:

\begin{eqnarray*}
\kappa (\psi_1,\ldots, \psi_L) &=& E[[2 I_p - E\{(f_t - g_s)^ {\otimes 2} | y_{t+1}, \eta_{s+1}\} - 2 E\{E^{\otimes 2}(f_t | \psi_1' f_t,\ldots, \psi_L' f_t) \}  \\
&& \hspace{0.35cm} + E[\{E(f_t | \psi_1' f_t,\ldots, \psi_L' f_t) -  E(g_s | \psi_1' g_s,\ldots, \psi_L' g_s)\}^ {\otimes 2} | y_{t+1}, \eta_{s+1}]]^ {\otimes 2}
\end{eqnarray*}
over all the semi-orthogonal matrices $(\psi_1,\ldots, \psi_L)$, where $v^{\otimes 2}$ denotes $v v'$ for any real vector $v$ and $E(f_t | \psi_1' f_t,\ldots, \psi_L' f_t)$ is modeled parametrically as if (B1') held for $(\psi_1,\ldots, \psi_L)$. Using the estimated factors $\widehat f_t$ and $\widehat g_s$ and the slicing strategy, we can similarly construct $\widehat \kappa (\cdot)$.

Under fairly general assumptions \citep{dong2010},  there exists the unique minimizer of $\kappa (\cdot)$ up to orthogonal column transformations, which spans the central subspace ${\mathcal S}_{y | f}$; we omit these assumptions here. Intuitively, a minimizer of $\widehat \kappa (\cdot)$ spans a consistent estimator of ${\mathcal S}_{y | f}$.

\begin{theorem}\label{thm:relaxB1}
Let $(\widehat \phi_1,\ldots, \widehat \phi_L)$ denote any minimizer of $\widehat \kappa(\psi_1,\ldots, \psi_L)$. Under Assumptions 1 -- 3, condition (B1'), and the constant variance condition (B2), we have
\begin{equation*}
\| (\widehat \phi_1,\ldots,\widehat \phi_L) (\widehat \phi_1,\ldots,\widehat \phi_L)' -  (\phi_1,\ldots, \phi_L)  ( \phi_1,\ldots, \phi_L)'\|_F = \bop (p^{-1/2} + T^{-1/2}).
\end{equation*}
\end{theorem}

By Theorem \ref{cor:dr-consist} and Theorem \ref{thm:relaxB1}, we can apply the proposed forecasting method or its generalization without concerning the linearity condition (B1), for both fixed and diverging $K$. For example, we now allow the predictor $x_t$, as well as the factors $f_t$, to contain discrete components.

\section{Determining $K$ and $L$}

We now discuss how to determine the number of factors $K$ and the dimension $L$ of the central subspace $\cstrue$. The problem is commonly called order determination in the literature of dimension reduction \citep{luo2016ladle}.

In the literature, various order-determination methods have been proposed to estimate $K$, including \cite{bai-2002, bai-ng-2008, onatski-2010, ahn-horenstein-2013}, \cite{Ludvigson-Ng-2009}, \cite{jurado-2015}. Recently, \cite{li-etal-2013} extended \citeauthor{bai-2002}'s approach to the case of diverging $K$, and estimated $K$ by
\begin{equation*}
\widehat K = \arg\min_{0\leq k \leq K_{\max}}\log(p\inv T\inv \|X - T^{-1}X\hF_k\hF_k'\|^2_F) + k \cdot q(p,T),
\end{equation*}
where $K_{\max}$ is a prescribed upper bound that possibly increases with $p$ and $T$, and {$\hF_k$ denotes the solution to (\ref{eq:est B and F}) with $k$ being the working number of factors}. $q(p,T)$ is a penalty function such that $q(p,T)=o(1)$ and {$(K_{\max}^6/p + K_{\max}^4/T)^{-1} q(p,T) \rightarrow \infty$}. We adopt \citeauthor{li-etal-2013}'s approach, and follow their suggestion to take
\begin{equation*}
q(p,T) = (p+T)(pT)\inv \log\{pT (p+T)\inv\}.
\end{equation*}

To estimate the dimension $L$ of the central subspace $\cstrue$, multiple methods have been proposed, including the sequential tests \citep{li-1991,li-wang-2007}, the bootstrap procedure \citep{ye2003using}, the cross-validation method \citep{xia-etal-2002,wang2008sliced}, the BIC-type procedure \citep{zhu-etal-2006}, and the ladle estimator \citep{luo2016ladle}, among which we adopt the BIC-type procedure and extend it to the high-dimensional case.
%For a sequence of kernel matrices $\{\bM_{T} \in \mathbb{R}^{K \times K}: T \in \mathbb{N}\}$ and its sample estimator $\{\widehat \bM_{T}: T \in \mathbb{N}\}$, denote the eigenvalues of each $\widehat \bM_{T}$
For {a positive semi-definite matrix parameter $M$ who columns span $\cstrue$} and its sample estimator $\widehat M$, let $\{\lambda_1,\ldots, \lambda_K\}$ and $\{\widehat \lambda_1,\ldots, \widehat \lambda_K\}$ be their eigenvalues in the descending order, respectively. {By definition, $\lambda_L$ must be positive.} We introduce a constant $c \in (0, 1)$ {and set $K_c$, the nearest integer to $cK$, as an upper bound of $L$. This is reasonable because $L$ is fixed and usually small in practice. }
We modify the objective function in \cite{zhu-etal-2006} to $G: \{1,\ldots, K_c\} \rightarrow \mathbb{R}$ with
\begin{equation}\label{eq:bic}
G(l) = (T/2) \tsum_{i=1 + \min (\tau, l)}^{K_c} \{\log (\widehat \lambda_i + 1) - \widehat \lambda_i\} - C_T l (2K - l + 1) / 2,
\end{equation}
where $\tau$ is the number of positive $\widehat \lambda_i$'s. We then estimate $L$ as the maximizer $\widehat L$ of $G(\cdot)$. Due to the introduction of the non-trivial upper bound $K_c$, we do not need to impose additional constraints on $K$ or $\|\widehat M - M\|$ for the consistency of $\widehat L$. This improves the result in \cite{zhu-etal-2006}.
%As $\tau$ is easily seen to converge to $K$ in probability, we should not worry about it in implementation.
%The consistency of $\widehat L$ in the following improves the original result in \cite{zhu-etal-2006}, as we impose no additional constraints on $K$ or $\|\widehat \bM - \bM\|$. This is due to  of $L$.

\begin{theorem}\label{thm:order-determination}
Suppose $\|\widehat M - M\| = \sop (1)$. If $C_T$ satisfies
$C_T K T^{-1} \rightarrow 0$ and $\|\widehat M - M\|^2 = \sop (C_T K T^{-1})$,
then $\widehat L$ converges to $L$ in probability.
\end{theorem}

A candidate of $C_T$ is $K^{-1} T \|\widehat M - M\|$. Referring to Theorem \ref{thm:dr-consist}, if we apply the BIC-type procedure to directional regression, then we can choose $C_T$ to be $K^{1/2} p^{-1/2} T + T^{1/2}$. %To further polish the procedure, we can incorporate a multiplicative constant in $C_T$, and tune its value in a data-driven manner such as cross-validation.

\section{Simulation studies}

We now present a numerical example to illustrate the performance of the proposed forecasting method that uses directional regression in the sufficient dimension reduction stage. The data generating process is specified as the following:
\begin{equation*}
y_{t+1}=g(\phi_1'f_t,\phi_2'f_t) + \sigma\epsilon_{t+1}, \quad \text{and} \quad x_{it}= b_i'f_t + u_{it},
\end{equation*}
We fix $\phi_1 = (1,1,1,0'_{K-3})/\sqrt{3}, \phi_2=(1,0'_{K-3},1,3)/\sqrt{11}$. Following \cite{li-etal-2013}, we set the number of factors $K$ to increase with $p$ in the form of $K = [1.5\log (p)]$, where $[x]$ denotes the integer part of a real number $x$.
%where the number of factors $K$ is taken to be six and we fix $\phi_1 = (1,1,1,0,0,0)/\sqrt{3}, \phi_2=(1,0,0,0,1,3)/\sqrt{11}$.
The factor loadings $b_i$ are independently sampled from $U[-1,2]$. We generate the latent factors $f_{j,t}$ and the error terms $u_{it}$ from two AR(1) processes, $f_{j,t} = \alpha_jf_{j,t-1}+e_{jt}$ and $u_{it} = \rho_iu_{i,t-1}+\nu_{it}$, with $\alpha_j,\rho_i$ drawn from $U[0.2,0.8]$ and fixed during the simulation, and the noises $e_{jt},\nu_{it}$, are $N(0,1)$. We set $\epsilon_{t+1}\sim N(0,1)$ and $\sigma=0.2$.

We consider four different choices of the link function $g(\cdot)$,

\noindent
Model I: $y_{t+1} = 0.4(\phi_1'f_{t})^2 + 3\sin(\phi_2'f_{t}/4) + \sigma\epsilon_{t+1}$;

\noindent
Model \two: $y_{t+1} = 3\sin(\phi_1'f_{t}/4) + 3\sin(\phi_2'f_{t}/4) + \sigma\epsilon_{t+1}$;

\noindent
Model \three: $y_{t+1} = 0.4(\phi_1'f_{t})^2 + |\phi_2'f_{t}|^{1/2} + \sigma\epsilon_{t+1}$;

\noindent
Model \four: $y_{t+1} = (\phi_1'f_{t})(\phi_2'f_{t}+1) + \sigma\epsilon_{t+1}$.

\noindent
The proposed forecasting by directional regression (DR) is compared with the forecasting by sliced inverse regression (SIR) \citep{fxy-2015}, the linear PC-estimator (principal components), and the semi-parametric efficient estimator (SEE) proposed by \cite{mazhu2013}. Model I and \three \ includes at least one symmetric component, which cannot be estimated well by SIR. Model \two \ is favorable to SIR. Model \four \ contains the interaction component to examine the ability of each method in detecting such nonlinear effect.

To gauge the quality of the estimated directions, we adopt the squared multiple correlation coefficient $R^2(\hphi) = \max_{\phi\in \cstrue,\|\phi\|=1}(\phi'\hphi)^2$, where $\cstrue$ is spanned by $\phi_1$ and $\phi_2$. We ensure that the true factors and loadings meet the identifiability conditions by calculating $H$ such that $T^{-1}HF'FH'=I_K$ and $H^{-1}B'BH^{-1}$ is diagonal. The rotated central subspace is then understood as $H^{-1}\cstrue$, which is still denoted as $\cstrue$ (see \cite{fxy-2015}).

Table \ref{tab1} compares the estimation of SIR and DR in simulation studies, where the PC is omitted as it produces only one directional estimate. It is evident that DR has substantial improvement over SIR in model I, \three \  and \four, with higher $R^2(\hphi)$ and lower variance. This is not surprising as DR explores higher conditional moments and hence incorporates more information. SEE is slightly better than SIR in these cases, but it also fails to capture $\phi_2$ accurately, partially due to its semi-parametric nature which typically requires lengthy steps to converge. In model \two, SIR, DR and SEE yield comparable results. We also observe that DR has outstanding performance in small samples, which makes it favorable in practice.

\begin{table}[!htbp]
\footnotesize
\begin{center}
\begin{tabular}[h]{cc|cc |cc|cc}
\hline
\hline
\multicolumn{2}{c}{Model I}&  \multicolumn{2}{c}{SIR} &  \multicolumn{2}{c}{DR} & \multicolumn{2}{c}{SEE} \\ % & \multicolumn{2}{c}{CSS} \\
$p$ & $T$ & $R^2(\hphi_1)$ & $R^2(\hphi_2)$ & $R^2(\hphi_1)$ & $R^2(\hphi_2)$ & $R^2(\hphi_1)$ & $R^2(\hphi_2)$ \\ % & $R^2(\hphi_1)$ & $R^2(\hphi_2)$ \\
\hline
100 & 100 & 75.0(21.3) & 28.4(27.4) & 82.9(14.8) & 79.9(21.9) & 80.4(26.8) & 27.0(23.5) \\ %& 60.5(29.4) & 34.5(28.0)  \\
100 & 200 & 88.7(10.4) & 17.7(27.6) & 94.5(5.4)  & 91.5(8.5)  & 83.4(26.6) & 21.7(20.7) \\ %& 73.5(28.1) & 41.9(29.1) \\
100 & 500 & 95.9(3.6)  & 14.4(28.2) & 98.4(1.4)  & 96.0(3.4)  & 87.6(26.9) & 30.8(20.7) \\ %& 83.3(23.7) & 57.8(26.7) \\
200 & 100 & 63.2(24.5) & 26.6(24.8) & 74.6(20.3) & 67.9(24.4) & 40.9(23.4) & 13.0(18.5) \\ %& 56.4(27.6) & 32.3(23.5) \\
500 & 200 & 76.6(16.1) & 16.1(23.2) & 86.8(20.1) & 80.2(22.1) & 26.6(15.4) & 9.4(15.8)  \\ %& 64.0(27.2) & 35.1(22.2) \\
500 & 500 & 90.5(5.5)  & 9.2(22.4)  & 96.0(29.9) & 87.7(26.0) & 24.2(13.4) & 7.6(13.5)  \\ %& 76.2(24.2) & 52.0(24.3) \\
\hline
\hline
\multicolumn{2}{c}{Model \two}&  \multicolumn{2}{c}{SIR} &  \multicolumn{2}{c}{DR} & \multicolumn{2}{c}{SEE} \\ %&  \multicolumn{2}{c}{CSS} \\
\hline
100 & 100 & 95.8(3.5) & 21.0(25.7) & 95.8(3.5)  & 26.4(26.6) & 89.7(22.5) & 33.0(20.1) \\ %& 60.8(30.7) & 30.8(25.6) \\
100 & 200 & 97.8(1.8) & 32.4(27.7) & 97.9(1.8)  & 43.4(28.7) & 90.4(15.0) & 30.2(19.5) \\ %& 77.1(29.8) & 37.2(27.9) \\
100 & 500 & 99.1(0.7) & 63.8(27.0) & 99.1(0.7)  & 74.8(23.8) & 91.9(20.5) & 48.7(20.9) \\ %& 85.8(25.3) & 56.0(27.3) \\
200 & 100 & 94.6(3.6) & 17.6(22.4) & 94.2(10.6) & 21.4(23.4) & 81.6(26.8) & 21.2(18.7) \\ %& 60.5(27.7) & 27.2(22.5) \\
500 & 200 & 95.9(2.1) & 18.2(22.6) & 95.5(11.9) & 24.7(23.4) & 37.8(26.5) & 13.9(17.4) \\ %& 67.6(25.5) & 34.2(22.6) \\
500 & 500 & 98.4(0.9) & 41.1(25.6) & 97.9(15.2) & 48.3(26.3) & 30.7(24.9) & 13.1(17.1) \\ %& 79.1(24.0) & 45.8(24.7) \\
\hline
\hline
\multicolumn{2}{c}{Model \three}&  \multicolumn{2}{c}{SIR} &  \multicolumn{2}{c}{DR} & \multicolumn{2}{c}{SEE} \\ %&  \multicolumn{2}{c}{CSS} \\
\hline
100 & 100 & 33.4(26.7) & 26.1(23.4 )& 83.0(19.7) & 47.6(28.2) & 40.1(30.7) & 29.9(18.4) \\ %& 67.1(28.8) & 50.2(31.8) \\
100 & 200 & 34.8(27.3) & 23.8(22.7) & 94.9(4.1)  & 83.2(22.9) & 68.4(35.1) & 20.2(18.1) \\ %& 87.1(26.8) & 85.0(29.5) \\
100 & 500 & 33.0(28.1) & 24.2(23.4) & 98.4(1.4)  & 97.6(2.1)  & 77.2(34.6) & 21.5(16.7) \\ %& 97.3(21.1) & 95.9(23.7) \\
200 & 100 & 29.5(25.9) & 19.8(20.4) & 75.0(23.3) & 36.5(25.7) & 37.9(26.8) & 12.9(17.9) \\ %& 67.1(28.2) & 44.5(30.2) \\
500 & 200 & 20.3(23.7) & 15.2(10.1) & 88.9(22.2) & 48.8(27.8) & 20.5(16.1) & 8.6(14.5)  \\ %& 75.6(27.1) & 58.6(30.2) \\
500 & 500 & 21.3(23.1) & 14.5(18.1) & 95.6(29.6) & 92.9(28.0) & 14.0(13.5) & 6.6(13.7)  \\ %& 91.5(24.5) & 85.9(29.7) \\
\hline
\hline
\multicolumn{2}{c}{Model \four}&  \multicolumn{2}{c}{SIR} &  \multicolumn{2}{c}{DR} & \multicolumn{2}{c}{SEE} \\ %&  \multicolumn{2}{c}{CSS} \\
\hline
100 & 100 & 61.8(29.1) & 31.3(26.0) & 85.6(14.2) & 79.1(23.5) & 64.4(27.8) & 43.9(18.4) \\ %& 85.5(21.5) & 85.5(20.0) \\
100 & 200 & 75.1(26.4) & 41.6(27.9) & 94.5(4.9)  & 93.5(5.2)  & 71.7(34.1) & 51.1(19.4) \\ %& 94.2(21.3) & 93.8(22.7) \\
100 & 500 & 89.4(15.0) & 67.8(27.4) & 98.1(1.7)  & 97.7(1.9)  & 88.2(37.7) & 66.6(17.0) \\ %& 97.2(20.8) & 97.0(22.4) \\
200 & 100 & 51.9(28.6) & 29.0(24.8) & 79.6(19.7) & 71.0(24.3) & 41.5(25.9) & 12.2(18.2) \\ %& 81.5(22.7) & 82.0(24.7) \\
500 & 200 & 59.4(27.9) & 30.2(24.4) & 87.5(21.7) & 86.2(20.2) & 19.5(15.4) & 7.4(13.5)  \\ %& 88.0(19.4) & 87.5(20.1) \\
500 & 500 & 83.3(17.8) & 54.9(26.9) & 95.1(28.3) & 94.6(26.4) & 10.2(13.3) & 4.8(13.1)  \\ %& 94.6(15.6) & 93.9(15.5) \\
\hline
\hline
\end{tabular}
%\vspace{0.3cm}
\caption{Performance of estimated $\hphi$ using median $R^2(\hphi)$ (\%) with standard deviations in parentheses over 1000 replications.
\label{tab1}}
\end{center}
\end{table}

\begin{table}[!ht]
\footnotesize
\begin{center}
\begin{tabular}{cc|cccc |cccc}
\hline
\hline
\multicolumn{2}{c}{}&  \multicolumn{4}{c}{Model I} &  \multicolumn{4}{c}{Model \two}\\
$p$ & $T$ & SIR & DR & PC & SEE & SIR & DR & PC & SEE \\
\hline
100 & 100 & -11.7 & 28.8 & -0.4 &  1.4 & 94.6 & 94.8 & 93.3 & 78.0 \\
100 & 200 &  -3.9 & 72.1 & 18.0 &  9.9 & 95.7 & 95.8 & 94.6 & 79.4 \\
100 & 500 &   0.4 & 92.2 & 27.4 & 11.4 & 96.1 & 96.2 & 94.9 & 79.3 \\
200 & 100 & -11.4 & 18.6 & -6.9 & -4.0 & 95.3 & 95.6 & 94.2 & 61.8 \\
500 & 200 &  -5.3 & 57.5 & -1.1 & -0.7 & 96.2 & 96.5 & 94.8 & 45.9 \\
500 & 500 &  -0.9 & 91.4 & 13.8 &  1.7 & 97.1 & 97.1 & 95.8 & 45.4 \\
\hline
\hline
\multicolumn{2}{c}{}&  \multicolumn{4}{c}{Model \three} &  \multicolumn{4}{c}{Model \four} \\
$p$ & $T$ & SIR & DR & PC & SEE & SIR & DR & PC & SEE \\
\hline
100 & 100 & -9.4 & 34.8 & 17.8 & 17.2 & -0.2 & 23.6 & 21.2 & 18.4 \\
100 & 200 &  1.0 & 77.1 & 30.8 & 22.7 & 13.5 & 53.7 & 35.8 & 28.4 \\
100 & 500 &  5.2 & 90.5 & 38.0 & 25.5 & 29.6 & 57.3 & 43.2 & 30.8 \\
200 & 100 & -9.7 & 21.5 &  3.8 &  6.3 & -2.3 & 16.9 &  6.8 &  7.6 \\
500 & 200 & -4.4 & 62.5 &  6.6 &  2.6 &  5.6 & 46.0 &  9.7 &  5.3 \\
500 & 500 & -1.3 & 89.5 & 19.1 &  5.2 & 22.4 & 58.3 & 21.6 & 48.5 \\
\hline
\hline
\end{tabular}
\caption{Comparison of out-of-sample median $R^2$
in percentage (\%) over 1000 replications. \label{tab2} }
\end{center}
\end{table}

We next investigate the predictive power of DR through the out-of-sample $R^2$, i.e.,
\begin{equation*}
R^2 = 1- {\tsum^{T+n_T}_{t=T+1}(y_t-\hat{y}_t)^2}/{\tsum^{T+n_T}_{t=T+1}(y_t-\bar{y}_t)^2},
\end{equation*}
where we use a fixed length $n_T=50$ of testing samples to evaluate the out-of-sample performance. $\hat{y}_t$ is the predicted value using all information prior to $t$. The fitting is done by building an additive model in Step 3 of the proposed estimator. In the case of PC-estimator, $\widehat{K}$ smooth functions are constructed for the estimated factors. In contrast, only $\widehat{L}$ smooth functions are applied in the cases of SIR, DR and SEE.
$\widehat{K}$ and $\widehat{L}$ are obtained using the procedures introduced in Section 5. It is clear from Table \ref{tab2} that DR enjoys great performance in almost all the cases. Similar to DR, SEE is better than SIR as it explores structural dimension more thoroughly with different forms of the target. But SEE is often limited to a large sample size to produce accurate estimation. The PC-estimator is more robust in the presence of symmetric components, but fails to capture the interaction effect in general. To investigate the accuracy of $\widehat K$ and $\widehat L$ used above, which are obtained from Section 5, we carry out simulations to investigate the accuracy of the estimation procedures, and examine the sensitivity of forecasting performance with respect to $\widehat{K}$ and $\widehat{L}$. In addition, we conduct experiments to show the effectiveness of the proposed method when the linearity condition is violated for factors $f_t$. Due to space limit, these numerical results are presented in the supplementary materials.

\section{Macro Index Forecast}

We now analyze how the diffusion indices constructed by the proposed DR impact real-data forecasts. We use a monthly macro dataset consisting of 134 macroeconomic time series recently composed by \cite{mccracken2016fred}, which are classified into 8 groups : (1) output and income, (2) labor market, (3) housing, (4) consumption, orders and inventories, (5) money and credit, (6) bond and exchange rates, (7) prices, and (8) stock
market. The dataset spans from 1959:01 to 2016:01. For a given target time series, we model the multi-step-ahead variable as:
$$
y_{t+h}^h=g(\phi_1'f_t,...,\phi_L'f_t)+\epsilon_{t+h}^h,
$$
where $y_{t+h}^h = h^{-1}\sum_{i=1}^h y_{t+i}$ is the variable to forecast, as in \cite{stock-watson-2002}.%, and $L$ is taken to be 1 or 2.

We follow \cite{mccracken2016fred} to preprocess the data. We also employ the Ljung-Box test with various lags to test for uncorrelatedness in residuals, which suggests the appropriateness to use our proposed methods. Forecasts of $y_{t+h}^h$ are constructed based on a moving window with fixed length ($T=120$) to account for timeliness. For each fixed window, the factors in the forecasting equation are estimated by the method
of principal components using all time series except the target. As noted by \cite{mccracken2016fred}, 8 factors have good explanatory power in various cases, so we set $K=8$ throughout the exercise. For each method $M$, we compare out-of-sample forecasting performances using the relative MSE (RMSE) to the PC method,
\begin{equation*}
\mathrm{RMSE}(M) = \mathrm{MSE}(M) \, / \, \mathrm{MSE}(\mathrm{PC}), \mbox{ where } \ \mathrm{MSE}(M)  = m\inv \tsum_{t=T+1}^{T+m}(y_t - \hat{y}_t)^2,
\end{equation*}
which we evaluate on the last $m=240$ months (20 years). The methods we consider here include SIR($i$), DR($i$) ($i=1,2$), where SIR($i$) denotes sufficient forecasting with $L=i$, and similarly for DR. Both methods use an additive model in specifying the forecasting equation. We also impose an additive model to the estimated factors, denoted by NL-PC, to see how much we can leverage on the nonlinearity without projecting principal components.

\begin{table}[!ht]
\footnotesize%
\begin{center}
\begin{tabular}{l|ccccc}
\hline
\hline
Group ($h=1$) &  SIR(1) & SIR(2) & DR(1) & DR(2) & NL-PC\\
\hline
Output \& Income & 1.03/1.61/0.96 & 1.02/1.13/0.94 & 0.99/1.19/0.92 & 1.02/1.14/0.90 & 1.21/1.38/1.05\\
Consumption & 1.00/2.10/0.80 & 0.95/1.05/0.74 & 0.92/1.02/0.86 & 1.00/1.05/0.81 & 1.16/1.44/1.04\\
Labor market & 1.02/2.27/0.71 & 1.00/1.21/0.42 & 0.97/1.13/0.52 & 0.98/1.16/0.42 & 1.21/1.53/0.46\\
Housing & 1.04/1.32/0.64 & 0.92/1.08/0.52 & 0.83/1.04/0.50 & 0.79/0.94/0.44 & 0.83/0.97/0.49\\
Money \& Credit & 0.94/1.04/0.86 & 0.97/1.05/0.90 & 0.96/1.10/0.86 & 1.04/1.24/0.92 & 1.14/1.41/1.07\\
Stock market & 0.99/1.39/0.90 & 1.02/1.12/0.83 & 0.92/1.08/0.88 & 1.04/1.07/0.91 & 1.36/1.39/1.14\\
Interest rates & 1.04/1.79/0.79 & 0.93/1.17/0.61 & 0.90/1.04/0.59 & 0.92/1.15/0.62 & 1.12/1.32/0.73\\
Prices & 0.97/1.42/0.80 & 0.99/1.05/0.83 & 0.95/1.12/0.81 & 0.97/1.12/0.88 & 1.12/1.47/0.92\\
\hline
\hline
Group ($h=6$)&  SIR(1) & SIR(2) & DR(1) & DR(2) & NL-PC\\
\hline
Output \& Income & 1.07/1.47/0.93 & 0.97/1.23/0.81 & 0.99/1.18/0.89 & 1.05/1.27/0.95 & 1.28/1.52/0.97\\
Consumption & 1.16/1.73/0.90 & 0.90/1.12/0.67 & 0.94/1.16/0.71 & 1.03/1.14/0.73 & 1.28/1.66/0.77\\
Labor market & 1.15/2.02/0.68 & 0.89/1.22/0.39 & 0.90/1.26/0.48 & 0.98/1.39/0.43 & 1.24/1.42/0.45\\
Housing & 0.96/1.29/0.66 & 0.85/0.95/0.51 & 0.73/0.89/0.50 & 0.69/0.86/0.47 & 0.78/1.02/0.55\\
Money \& Credit & 0.95/3.51/0.76 & 1.01/3.65/0.83 & 0.99/1.52/0.76 & 1.02/1.74/0.78 & 1.23/2.90/0.92\\
Stock market & 0.91/1.20/0.83 & 0.94/1.05/0.89 & 0.89/1.08/0.84 & 1.00/1.03/0.94 & 1.23/1.27/0.83\\
Interest rates & 1.01/1.61/0.75 & 0.90/1.12/0.64 & 0.84/1.13/0.50 & 0.88/1.18/0.58 & 1.11/1.46/0.70\\
Prices & 1.16/1.37/0.51 & 1.03/1.12/0.82 & 1.11/1.37/0.94 & 1.14/1.36/0.95 & 1.17/1.35/1.11\\
\hline
\hline
Group ($h=12$) &  SIR(1) & SIR(2) & DR(1) & DR(2) & NL-PC\\
\hline
Output \& Income & 1.24/1.67/0.79 & 1.01/1.45/0.76 & 0.99/1.22/0.76 & 1.01/1.36/0.86 & 1.17/1.34/0.92\\
Consumption & 1.27/1.60/0.83 & 1.08/1.44/0.62 & 1.09/1.32/0.65 & 1.06/1.38/0.66 & 1.16/1.38/0.87\\
Labor market & 1.07/1.76/0.67 & 0.83/1.40/0.41 & 0.91/1.44/0.54 & 0.89/1.41/0.46 & 1.13/1.39/0.56\\
Housing & 0.85/1.35/0.59 & 0.69/0.93/0.46 & 0.67/0.91/0.40 & 0.68/0.83/0.36 & 0.89/1.16/0.54\\
Money \& Credit & 1.14/2.03/0.41 & 1.03/2.16/0.80 & 1.05/1.52/0.85 & 1.00/1.40/0.82 & 1.20/1.69/0.87\\
Stock market & 1.09/1.20/0.89 & 1.01/1.13/0.84 & 0.96/1.17/0.94 & 1.08/1.16/0.75 & 1.06/1.14/0.89\\
Interest rates & 1.00/1.31/0.75 & 0.82/1.22/0.59 & 0.80/1.27/0.53 & 0.85/1.18/0.51 & 1.07/1.62/0.70\\
Prices & 1.18/1.40/0.53 & 1.21/1.40/0.66 & 1.19/1.31/0.71 & 1.21/1.33/0.77 & 1.25/1.52/0.94\\
\hline
\hline
\end{tabular}
\footnotesize

 \caption{RMSE in Out of Sample Forecast (Median/Max/Min): out-of-sample RMSE relative to the linear diffusion index. In each group, the median, maximum and minimum of RMSE is reported. SIR($i$) denotes sufficient forecasting using $i$ indices, DR denotes sufficient directional forecasting, and NL-PC denotes a nonlinear additive model on all the estimated factors. \label{tab3}}
\end{center}
\end{table}

\begin{figure}[!htbp]
\begin{center}
\includegraphics[scale=0.36]{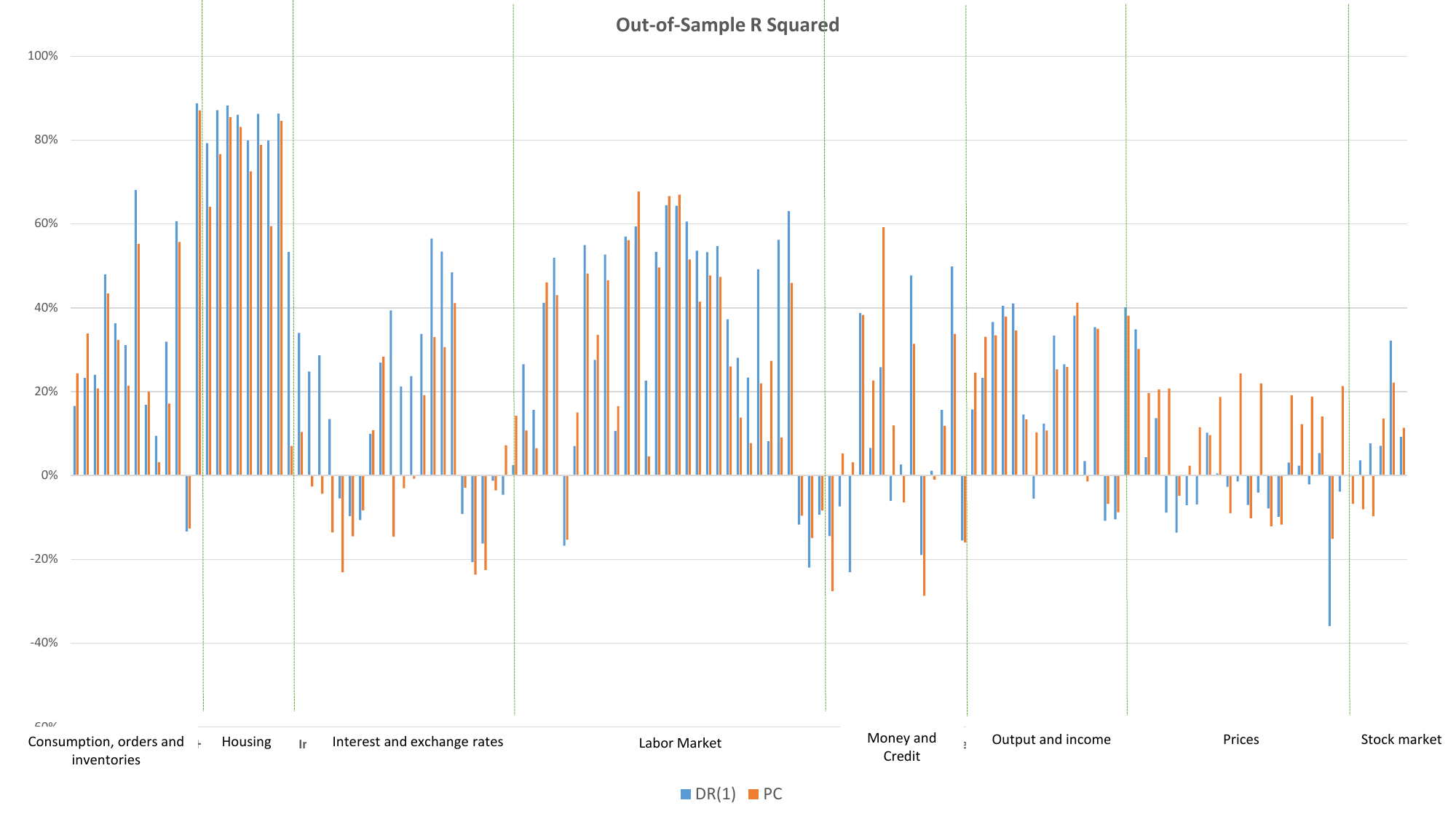}
\caption[Out-of-Sample $R^2$]{\small 6-month ahead forecasting out-of-sample $R^2$ for the 134 macroeconomic series organized into eight groups.%, using DR(1) and PC methods.
}
\end{center}
\label{dr1}
\end{figure}

\medskip

We report results in Table \ref{tab3} for $h=1,6,12$, on the maximum, minimum and median of RMSE in each broad sector. Several features are noteworthy. First, a nonlinear additive model {built on estimated factors} does not buy us more predictive power, except in the housing sector, where most of the nonlinear methods improve prediction accuracy. Second, the one-step-ahead out-of-sample forecast favors DR(1), as we observe the median RMSEs are uniformly less than 1 and some of the reductions in RMSE are substantial. Moving from short horizon to long horizon changes predictability of the targets, but DR(1) manages to improve the forecast over the PC method in many instances. Finally, as an illustration, we plot the out-of-sample $R^2$ for the 6-month-ahead forecast using DR(1) and PC. Notably, macro time series in housing and labor market sectors have higher predictability than in rates and stock market sectors.

\end{document}